%%%%%%%%%%%%%%%%%%%%   Geometry and Topology: 1998-4.tex  %%%%%%%%%%%%%
%%%%        
%%%%  Group negative curvature for 3-manifolds with genuine laminations
%%%%             
%%%%                  David Gabai, William H. Kazez
%%%%  
%%%%             Published in Volume 2(1998) pages 65--77 
%%%%
%%%%                   Publication date 11 May 98 
%%%%
%%%%                   This is a plain TeX file
%%%%
%%%%
%%%%%%%%%%%%%%%%%%              %%%%%%%%%%%%%%%%%%
%
%%%%%%%%%%%%%%%%%%%%%%   amstex.pat   %%%%%%%%%%%%%%%%%%%
%%
%%    Patch to adapt an old-fashioned amstex file to fit
%%    GT style.  Insert the patch at the start in place of
%%    \documentstyle{amsppt}.  
%%
%%    Use plain tex (NOT amstex) to process the file.
%%    ===============================================
%%
%%    Note that \thanks and \dedicatory are disabled.  So be sure 
%%    to move any acknowledgements "manually" to the end of the 
%%    introduction.  Other than dealing with the references
%%    (see macros below) this should be all the manual work
%%    needed.
%%   
%%    One other hack needed --- search for \@ and replace by @
%%
%%                                        Colin Rourke 28-8-97
%%                                             updated 11-4-98
\input gtmacros
\input amsnames
\input amstex
\input rlepsf         %%% uncomment if needed and add further 
\let\cal\Cal          % amstex disables \cal 
\catcode`\@=12        % and leaves @ \active 
%
%                    Publication information (uncomment)
%
\input gtoutput
\volumenumber{2}\papernumber{4}\volumeyear{1998}
\pagenumbers{65}{77}\published{11 May 1998}
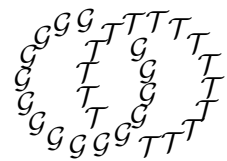
\shorttitle{Group negative curvature for laminar 3-manifolds}  
\proposed{Jean-Pierre Otal}\seconded{Robion Kirby, Michael Freedman}
\received{5 August 1997}\revised{9 May 1998}
%\accepted{24 February 1998}
%
%%
%%
%%     Adapt basic layout macros:
%%
\let\\\par
\def\topmatter{\relax}
\def\endtopmatter{\maketitlepage}
\let\gttitle\title
\def\title#1\endtitle{\gttitle{#1}}
\let\gtauthor\author
\def\author#1\endauthor{\gtauthor{#1}}
\let\gtaddress\address
\def\address#1\endaddress{\gtaddress{#1}}
\let\gtemail\email
\def\email#1\endemail{\gtemail{#1}}
\def\subjclass#1\endsubjclass{\primaryclass{#1}}
\let\gtkeywords\keywords
\def\keywords#1\endkeywords{\gtkeywords{#1}}
\def\heading#1\endheading{{\def\S##1{\relax}\def\\{\relax\ignorespaces}
    \section{#1}}}
\def\head#1\endhead{\heading#1\endheading}

\def\subhead#1\endsubhead{\sh{#1}}
\def\subsubhead#1\endsubsubhead{\sh{#1}}
\def\specialhead#1\endspecialhead{\sh{#1}}
\def\demo#1{\rk{#1}\ignorespaces}
\def\enddemo{\ppar}

\def\qed{\ifmmode\quad\sq\else\hbox{}\hfill$\sq$\par\goodbreak\rm\fi}  
\def\proclaim#1{\rk{#1}\sl\ignorespaces}
\def\endproclaim{\rm}
\def\cite#1{[#1]}
\newcount\itemnumber
\def\roster{\items\itemnumber=1}
\def\endroster{\enditems}
\let\itemold\item
\def\item{\itemold{{\rm(\number\itemnumber)}}%
\global\advance\itemnumber by 1\ignorespaces}
\def\S{section~\ignorespaces}  %%  expand \S to "section"
\def\date#1\enddate{\relax}
\def\thanks#1\endthanks{\relax}   %%%  Move acknowledgements "manually"
\def\dedicatory#1\enddedicatory{\relax}  %%% to the end of the intro.
  % in some versions of amstex but not all
%%
%%
%%   Adapt the amstex reference list macros 
%%   (some spacing may need "manual" adjustment) :
%%
%%
\def\Refs{\ppar{\large\bf References}\ppar\bgroup\leftskip=25pt
\frenchspacing\parskip=3pt plus2pt\small}       
\def\endRefs{\egroup}
\def\widestnumber#1#2{\relax}
\def\endrefitem{}
\def\refdef#1#2#3{\def#1{\endrefitem#2\def\endrefitem{#3}}}
\def\ref{\par}
\def\endref{\endrefitem\par\def\endrefitem{}}
\refdef\key{\noindent\llap\bgroup[}{]\ \ \egroup}
\refdef\no{\noindent\llap\bgroup[}{]\ \ \egroup}
\refdef\by{\bf}{\rm, }
\refdef\manyby{\bf}{\rm, }
\refdef\paper{\it}{\rm, }
\refdef\book{\it}{\rm, }
\refdef\jour{}{ }
\refdef\vol{}{ }
\refdef\yr{$(}{)$ }
\refdef\ed{(}{ Ed.) }
\refdef\publ{}{ }
\refdef\inbook{from: ``}{'', }
\refdef\pages{}{ }
\refdef\page{}{ }
\refdef\paperinfo{}{ }
\refdef\bysame{\hbox to 3 em{\hrulefill}\thinspace,}{ }
\refdef\toappear{(to appear)}{ }
\refdef\issue{no.\ }{ }
%%
%%      Macros to change refs to numbers
%%
%%      To use these macros uncomment the macros below :
%%
\newcount\refnumber\refnumber=1
\def\refkey#1{\expandafter\xdef\csname cite#1\endcsname{\number\refnumber}%
\global\advance\refnumber by 1}
\def\cite#1{[\csname cite#1\endcsname]}
\def\Cite#1{\csname cite#1\endcsname}  %% unbracketed \cite 
\def\key#1{\noindent\llap{[\csname cite#1\endcsname]\ \ }}
%%
%%      Next edit the reference list so that all keys
%%      are enclosed in { }'s.  Then copy
%%      the list below and edit and cut the lines down to
%%      produce a \refkey list, of which the following 
%%      are sample lines:
%%
%%      \refkey {G2}
%%      \refkey {G3}
%%
\refkey {BF}
\refkey {Br1}
\refkey {Br2}
\refkey {Cl}
\refkey {DR}
\refkey {G1}
\refkey {G2}
\refkey {G3}
\refkey {GK1}
\refkey {GK2}
\refkey {GO}
\refkey {Gr}
\refkey {Ha}
\refkey {J}
\refkey {JS}
\refkey {M}
\refkey {Na1}
\refkey {Th1}
\refkey {Th2}
\refkey {Wu}
%%%%
%%       End of patch
%%

\define\cirE{\ \raise4pt\hbox{$^\circ$}\kern-8pt E}

\define\finv{f^{-1}}

\define\BR{\Bbb R}

\define\F {\Cal F}
\define\T {\Cal T}

\TagsOnRight

\topmatter

\title Group negative curvature for 3--manifolds\\with genuine laminations
\endtitle

\author David Gabai\\William H. Kazez
\endauthor 

\address California Institute of Technology\\Pasadena, CA
91125-0001 USA\endaddress

\email gabai@cco.caltech.edu\\will@math.uga.edu\endemail

\secondaddress{University of Georgia\\Athens, GA 30602, USA}

\abstract

We show that if a closed atoroidal 3--manifold $M$ contains a genuine
lamination, then it is group negatively curved in the sense of Gromov.
Specifically, we exploit the structure of the non-product complementary
regions of the genuine lamination and then apply the first author's
Ubiquity Theorem to show that $M$ satisfies a linear isoperimetric inequality.

\endabstract

\asciiabstract{We show that if a closed atoroidal 3-manifold M 
contains a genuine
lamination, then it is group negatively curved in the sense of Gromov.
Specifically, we exploit the structure of the non-product complementary
regions of the genuine lamination and then apply the first author's
Ubiquity Theorem to show that M satisfies a linear isoperimetric 
inequality.}

\gtkeywords{Lamination, essential lamination, genuine lamination, group
negatively curved, word hyperbolic}

\primaryclass{57M50}
\secondaryclass{57R30, 57M07, 20F34, 20F32, 57M30}

\endtopmatter

\document

	\heading{\S 0 \\ Introduction}\endheading

In 1985 the notion of manifold with essential lamination was introduced
\cite{GO} to simultaneously generalize the notion of Haken manifold or
manifold with Reebless foliation.   A longstanding goal has been to generalize
to laminar manifolds  various results established for Haken manifolds.  See
\cite{GO}; \cite{Cl}, \cite{Br1}; \cite{Br2}, \cite{G2}; for some results in
this direction. In particular, a major goal is to show that atoroidal laminar
manifolds have metrics of constant negative curvature, thereby generalizing
Thurston's hyperbolization theorem \cite{Th1}. Our main result shows that at
the level of the fundamental group, an atoroidal 3--manifold is negatively
curved (ie Gromov hyperbolic) provided that it has a genuine essential
lamination, ie a lamination which is not a Reebless foliation split open
along leaves.  

\proclaim {Theorem~0.1}  If $M$ is a closed atoroidal genuinely laminar
3--manifold,  then $\pi_1(M)$ is word hyperbolic.\endproclaim

\demo {Remarks 0.2} i)\stdspace  Recall that $M$ is {\sl atoroidal} if it contains
no immersed $\pi_1$--injective torus, equivalently $\pi_1(M)$ contains no
rank--2 free abelian subgroup.  

ii)\stdspace    Manifolds with essential laminations are far more plentiful then
manifolds with incompressible surfaces, eg see \cite{Na1},
\cite{DR}, \cite{Wu}, \cite{GO} and  \S 1 of the survey \cite{G1}.  Gabai and
Mosher \cite{M} showed that if $k$ is a hyperbolic knot in $S^3$, then off of a
line in Dehn surgery space, all manifolds obtained by surgery on $k$ have an
essential lamination.  As of this writing there is no known example of an
irreducible, atoroidal, 3--manifold with infinite $\pi_1$ which does not have
an essential lamination.  Although there do exist explicit examples not known
to have essential laminations.  See \S 1 of \cite{G1}.

 iii)\stdspace A {\sl genuine lamination} \cite{GK1} is an essential
lamination $\lambda$ which cannot be {\sl trivially} extended to a
foliation.  This means that some closed complementary region of
$M-\lambda$ is not an $I$--bundle.  Manifolds with genuine laminations
include manifolds with pseudo-Anosov flows, and by Gabai and Mosher
\cite{M} for each hyperbolic knot $k \in S^3$ off at most two lines in
its Dehn surgery space, all manifolds obtained by surgery.

Haken manifolds are compact orientable irreducible 3--manifolds that
contain an embedded incompressible surface.  Such a surface $S$ is
always an essential lamination and is a genuine lamination except in
the case that its closed complementary region is an $I$--bundle.  In
this case $S$ lifts to a fibre of a fibration in a 2--fold covering of
the manifold. If M is atoroidal, then by
\cite{Th2}, this covering has a pseudo-Anosov flow and hence a genuine
lamination.  Thus a closed atoroidal Haken manifold is 2--fold covered by a
manifold with a genuine essential lamination.

iv)\stdspace   It is the non $I$--bundle complementary region that allows one to
get a grip on the manifold and thereby extend to manifolds with genuine
laminations various known properties of Haken manifolds, see
\cite{GK1}, \cite{GK2}.  We  remarked in \cite{GK1} that the non
$I$--bundle part makes the manifold ``reek of negative curvature".  The
goal of this paper is to make that remark more precise.

v)\stdspace   Thurston's hyperbolization theorem asserts that atoroidal Haken
manifolds have hyperbolic structures.  Subsequently Bestvina, Feighn \cite{BF}
gave an elementary argument establishing group negative curvature for such
manifolds.  \enddemo

\demo{Idea of the Proof of Theorem~1}  Suppose that the genuine lamination 
$\lambda$ is obtained by splitting open a singular foliation $\mu$ along the
singular leaves where each singular leaf is a trigon $\times S^1$.  Let $C$ be
the singular locus of $\mu$.  The Ubiquity Theorem of
\cite{G3} asserts that there exists constants $K,L$ such that if $D$ is a least
area disk with $\text{length}(\partial D)/\text{area}(D)<L$, then $|E\cap
C|/\text{area}(D)>K$ where $E$ is any disk which spans $\partial D$.  For
reasonable disks, this means that after homotopy of $D$ rel $\partial D$, the
induced singular foliation on $D$ has at least $K \text{area}(D)$ 3--prong
complementary regions and $\mu|D$ has no circle leaves.  The Poincar\'e--Hopf
index formula then implies that $\partial D$ is tangent  to $\mu$
in at least $K \text{area}(D)/3$ different spots.  This in turn implies
that $\text{length}(\partial D)>K^\prime K \text{area}(D)/3$ and hence
$\text{length}(\partial D)/\text{area}(D)>K^\prime K/3$, for some constant
$K^\prime$ which is defined independently of reasonable $D$. 
For example, $K^\prime$ can be taken to be the minimal distance
between tangencies along $\partial D$.  In summary, $\text{area}(D)$ is
approximated by the geometric intersection number of $\partial D$ and $C$. 
The latter is approximately the number of tangencies of $\partial D$ with
$\mu$. This in turn gives a lower bound for $\text{length}(\partial D)$. 
Therefore, $M$ satisfies a linear isoperimetric inequality, and hence
$\pi_1(M)$ is negatively curved by Gromov \cite{Gr}.

The actual proof is not much different.  Since $\lambda$ is genuine,
there exists a finite, non-empty collection of characteristic annuli
embedded in $M-\lambda$, which separates off the $I$--bundle part of the
complement of $\lambda$ from the non $I$--bundle part (Lemma~1.3).  Taking $C$
to be the union of cores of these annuli, we obtain length and area
approximations in \S2 similar to those of the previous paragraph to conclude
that $M$ satisfies a linear isoperimetric inequality.

\enddemo

{\it  The authors would like to dedicate this paper to David
Epstein on the occasion of his 60th birthday.}

The first author was partially supported by
NSF Grant DMS-9505253 and the MSRI.

\heading{\S 1 \\Preliminaries}\endheading

\demo {Notation 1.1} Let $\raise4pt\hbox{$^\circ$}\kern-8pt E$ denote the
interior of $E$ and $|X|$ denote the number of connected components of $X$,
if $X$ is a space, or the number of elements of $X$ if it is just a set.
Let $\bar B$ denote the closure of $B$, and $I$ denote $[0,1]$.\enddemo 

In this section we remind the reader of several fundamental properties of
essential laminations.  

\demo {Definition 1.2}  A 2--dimensional {\sl lamination} $\lambda$ in a
3--manifold $M$ is a foliation of a closed subset of $M$.  More precisely
$M$ has charts of the form $\BR^2\times \BR$, such that
$\lambda\mid\BR^2\times \BR$ is the product lamination $\BR^2\times T$, 
where $T$ is a closed subset of $\BR$.  A component of $M-\lambda$ is called
a {\sl complementary region}.  A {\sl closed complementary region} is a
component $V$ of $M-\lambda$ metrically completed with respect to the
induced path metric. Topologically it is $V$ together with its boundary
leaves.  We will assume that the leaves of $\lambda$ are smoothly immersed
in $M$, although the transverse structure may only be $C^0$.

A lamination $\lambda$ in the closed orientable 3--manifold $M$ is {\sl
essential} \cite{GO} if there are no 2--sphere leaves, each leaf is
$\pi_1$--injective,  $M-\lambda$ is incompressible and each closed
complementary region is {\sl end-incompressible}.  The closed
complementary region $V$ is end-incompressible if for every proper map
$f\co D^2-x\to V$, $x\in \partial D^2$  such that $f(\partial D^2-x)\subset L$,
$L$ a leaf of $\lambda$, then there exists a proper map $g\co D^2-x\to L$ such
that $g\mid\partial D^2-x=f\mid\partial D^2-x$.

\enddemo

Let  $X$  be a codimension--$0$  submanifold of a 3--manifold  $V$.  The {\sl
horizontal boundary} of  $X$  is defined to be $X \cap
\partial V$ and is denoted $\partial_hX$.  The {\sl vertical boundary} of 
$X$  is defined to be the closure in $\partial X$  of  $\partial X -
\partial_hX$ and is denoted  $\partial_v X$.  Typically $V$ will be a closed
complementary region so that
$\partial_h X = X \cap \lambda$ and $\partial_vX$ will be a union of
annuli.  The pair $(X,\partial_h X)$  is called an
$I${\sl--bundle} if $X$  is the total space of an  $I$--bundle over a surface 
$S$  in such a way that  $\partial_h X$ consists of the boundary points of
the $I$--bundle fibres.

A {\sl product disk} is defined to be a proper embedding of  $(I
\times I, \partial I \times I, I \times \partial I)$  into 
$(X,\partial_vX, \partial_hX)$  such that each component of 
$\partial I \times I$ is mapped to an essential arc in 
$\partial_vX$.  A product disk is {\sl essential} if it is not parallel,
keeping  $\partial I \times I$ in  $\partial_vX$  and  $I
\times \partial I$  in  $\partial_hX$, to a disk in  $\partial_vX$.

\proclaim{Lemma~1.3} Let $V$ be the disjoint union of
the closed complementary regions of an essential lamination
$\lambda$.  There exists a unique (up to isotopy in $V$) finite collection
$\Cal A = A_1 \cup \dots \cup A_n$ of properly embedded annuli in $V$ such
that

\roster

\item $V = \Cal G \cup \Cal I$ where $\Cal G \cap \Cal I
=\partial_v \Cal G = \partial_v \Cal I = \Cal A$.

\item $(\Cal I, \partial_h \Cal I)$ is an $I$--bundle over a
possibly noncompact or disconnected surface.  No component of $(\Cal I,
\partial_h\Cal I)$ is an $I$--bundle over a compact surface with non-empty
boundary.

\item $(\Cal G, \partial_h \Cal G)$ is compact, has no
components homeomorphic to $(D^2 \times I, D^2 \times \partial I)$
and contains no essential product disks.
\endroster
\endproclaim

\demo{Proof} This follows almost immediately from \cite{J} or the Generalized
Splitting Theorem of \cite{JS}.   We shall only point out the minor
differences.

Their theorem is stated for compact pairs $(V, \partial_h V)$.  If
$N(B)$ is a fibred neighborhood of a branched surface which
carries $\lambda$, then $N(B)\cap V$ gives an $I$--bundle fibring of
the ends of $(V, \partial_h V)$.  Thus their result may be applied to
the finitely many components of the space $(V, \partial_h V)$ split along
$\partial_v (N(B)\cap V)$ which are compact.

Let $(\Cal I,\partial_h\Cal I)$ be those finitely many components of the
maximal $I$--bundle of $(V,\partial_h V)$ given by \cite{JS}, \cite{J} which
either contain an end of $V$ or are of the form $S\times I$
where $S$ is a closed surface.  Let $\Cal G=V- (\Cal I -\partial_v
\Cal I)$ and $\partial_h {\Cal G}={\Cal G}\cap \partial_h V$.  (3)
follows by the maximality of $\Cal I$.  By construction $\Cal I$ has
no $D^2 \times S^1$ components and consequently establishing uniqueness
is routine.\qed

\enddemo

Figures 1.1 and 1.2 show examples of the decomposition $V=\Cal G \cup \Cal I$.  Notice in Figure~1.1 that compact $I$--bundles are part of the guts $\Cal
G$.  In Figure~1.2 $\Cal I$ is maximal among $I$--bundles with no compact
components, but there are two different ways to add compact $I$--bundles to
$\Cal I$.
\medskip

\hskip .1truein
\relabelbox\small
\epsfxsize 5truein \epsfbox{Figure1.1}
\adjustrelabel <0pt,3pt> {G}{$\cal G$}
\relabel{I1}{$\cal I$}
\relabel{I2}{$\cal I$}
\relabel{A1}{$\cal A$}
\relabel{A2}{$\cal A$}
\relabel{I-bundle}{$I$--bundle\,\,over}
\relabel{compact}{compact\,\,subsurface}
\endrelabelbox
\vskip -.1truein
\centerline{\small Figure 1.1}
\medskip

\hskip .65truein
\relabelbox\small
\epsfxsize 4truein \epsfbox{Figure1.2}
\adjustrelabel <0pt,3pt> {I1}{$\cal I$}
\adjustrelabel <0pt,3pt> {I2}{$\cal I$}
\adjustrelabel <0pt,3pt> {A1}{$\cal A$}
\adjustrelabel <0pt,3pt> {A2}{$\cal A$}
\relabel{G1}{$\cal G$}
\relabel{G2}{$\cal G$}
\relabel{S1}{$\times S^1$}
\relabel{S2}{$\times S^1$}
\endrelabelbox

\centerline{\small Figure 1.2}
\medskip

\demo{Definition 1.4} Let $\Cal C$ denote the union of the cores of the
$A_i$'s.  $\Cal I$ together with its $I$--bundle structure is called the
{\sl interstitial bundle} of $\lambda$. $\Cal G$ is called the {\sl
guts} of $\lambda$.\enddemo

\demo {Definition 1.5}  A {\sl genuine lamination} is an essential lamination
in $M$ with non-empty guts.  In other words, a genuine lamination is an
essential lamination which is not just a Reebless foliation with some
leaves split open.  \enddemo

\proclaim {Lemma~1.6}  If $M$ has an essential lamination $\lambda_1$, then
$M$ has an essential lamination $\lambda$ which is nowhere dense and has no
isolated leaves.  If $\lambda_1$ is genuine, then $\lambda$ can be taken to
be genuine.\endproclaim  

\demo {Proof}  If $M-\lambda_1=\emptyset$, then by first Denjoy splitting
open a leaf we can assume that $\lambda_1$ is not all of $M$.  Let
$\lambda_2$ be a sublamination of $\lambda_1$ such that each leaf is dense in
$\lambda_2$. If $\lambda_2$ is a single compact leaf, then we may assume that
it is 2--sided by passing, if necessary, to a double cover.  Finally replace
it by $\lambda =\lambda_2 \times T$, where
$T$ is a Cantor set.  Use the fact that a sublamination of of an
essential lamination is essential by \cite{GO} to show that $\lambda$ is
essential.  Finally if $\lambda_1$ is genuine, then so is
$\lambda$.\qed\enddemo

\let\SS\stdspace
\demo {Remark 1.7}  i)\SS  The feature of no isolated leaves allows us to
assume that distinct $A_i$'s do not intersect.

ii)\SS  Let $\alpha$ be a properly embedded arc in a closed complementary
region $V$ of $\lambda$.  We say that $\alpha$ is {\sl efficient} if
$\alpha$ cannot be homotoped rel $\partial \alpha$ into
$\partial V$, via a homotopy supported in $V$. \enddemo

To establish Theorem~0.1 we need the following special case of Theorem~1
of \cite{GO}.  

\proclaim {Theorem~1.8} Let $\lambda$ be an essential lamination in $M$,
then 

i)\SS Interstitial fibres are efficient arcs.

ii)\SS  If I is an efficient arc, then $I$ cannot be homotoped (in $M$) rel
$\partial I$ to lie in a leaf.  

iii)\SS   If $C_i$ is a core of $A_i$, then $0\neq [C_i]\in \pi_1(M)$.

iv)\SS  If $\psi$ is a closed efficient transversal curve (ie $\psi$ is
transverse to $\lambda$ and the closure of each component of $\psi-\lambda$ is
an efficient arc), then $0\neq [\psi]\in \pi_1(M)$. \qed
\endproclaim

\heading{\S 2 \\ Proof of Theorem~0.1}\endheading

Let $\lambda$ be a genuine lamination in the closed 3--manifold $M$.  By
Lemma~1.6 we can assume that $\lambda$ is nowhere dense and that no leaf of
$\lambda$ is isolated.  As in \S 1 we let $\{A_i\}$ be the finite set of
characteristic annuli, $\Cal A=\cup A_i$ and let $\Cal C$ be the union of
the cores of the $A_i$'s.

Fix a Riemannian metric on $M$.  By Gromov \cite{Gr},  $\pi_1(M)$ is
word hyperbolic if there exists $L>0$ such that for each least area
mapping of a disk $f\co D\to M$, $\text{length}(\partial 
f)/\text{area}(f)>L$.

\proclaim{Lemma~2.1} Let $K >0$.  There exists $\epsilon >0$ with
the following property.  If the set of least area functions $f\co D\to
M$ satisfying

\roster
\item  $d(\partial f(\partial D), \Cal A)>\epsilon$,

\item  $\text{length}(\partial f|Y)>\epsilon$ for each component $Y$ of
$\finv(\Cal G)\cap
\partial  D$,

\item  $\partial f$ is transverse to $\lambda$ and

\item   $f(\partial D)$ intersects $\Cal I$ in interstitial
fibres

\endroster
all satisfy the isoperimetric inequality $\text{length}(\partial
f)/\text{area}(f)>K$, then $M$ is negatively curved.

\endproclaim

\demo{Proof} First observe that there exists $C_1$, and $\epsilon$ such
that if $\partial g_0\co S^1\to M$ is any smooth immersion, then  $\partial g_0$
is homotopic to $\partial g_1$ via $G\co S^1\times I\to M$ such that $\partial 
g_1$ satisfies the conclusions of (1)--(4), area$(G)\le C_1
\text {length}(\partial g_0)$ and  $\text{length}(\partial g_0)\ge \epsilon
\text {length}(\partial g_1)$.  (Hint: Think of  $\lambda$ as being
carried in an extremely thin fibred neighborhood of the smooth  branched
surface $B$, with the various components of $\Cal A$ lying in tiny
neighborhoods  of sparsely placed embedded curves in $B$.  The homotopy of
$\partial g_0$    to $\partial  g_1$ is more or less one which first makes
$\partial g_0$ piecewise geodesic, then pushes $\partial g_0$ off a small
neighborhood of $\Cal A$ and is transverse to $B$.)

Let $g_0\co D\to M$ be an arbitrary least area map.  Suppose that $\partial
g_1$ is as  above, and $g_1\co D\to M$ is any least area map which extends
$\partial g_1$.  Let $K_1 =\text{min}\{K, 1\}$.  Then
$$K<\text{length}(\partial g_1)/\text {area}(g_1)\le \frac{1}{\epsilon}
\text{length}(\partial g_0)/\text {area}(g_1)$$ 
and hence 
$$K_1< (\frac{1}{\epsilon}\text{length}(\partial g_0)+C_1\text {length} 
(\partial g_0))/(\text {area}(g_1)+C_1\text {length}(\partial g_0))$$
$$\le (\frac{1}{\epsilon}+C_1)\text {length}(\partial g_0)/\text
{area}(g_0)$$ which implies that  $\text{length}(\partial g_0)/\text
{area}(g_0)\ge K_1/(\frac{1}{\epsilon}+C_1)$. \qed

\enddemo

Let $E$ be any smooth branched immersed disk that spans $\partial D$
and minimizes geometric intersection number with $\Cal C$.  In what
follows (Claim~3) we establish a relationship between $|E\cap\Cal C|,\
|E\cap\Cal G|$ and $\text{length}(\partial D)$.

%%%% Editorial change : inserted \Cal twice before C in previous para --CPR

After a homotopy supported away from $\partial D$ we can assume that $E$ is
transverse to $\lambda \cup \Cal A$.  Since leaves of essential laminations
are $\pi_1$--injective and closed transverse efficient curves are
homotopically nontrivial, the induced lamination $\lambda|E$ is a lamination
by circles and properly embedded arcs.  Furthermore each circle is
homotopically trivial in its leaf.  Therefore the standard Reeb stability
argument implies the following lemma.

\proclaim{Packet Lemma~2.2} $\lambda|E$ is a union of product
laminations of the form $I\times T$ or $S^1\times T$ where $T$ is a Cantor
set in $I$ and $I\times I$ is an embedded square in $E$ with $\partial
E\cap I\times I=\partial I\times I$ or $S^1\times I$ is an embedded
annulus in $\cirE$.  Furthermore the various $I\times I$'s and the $S^1\times
I$'s are pairwise disjoint.\qed

\endproclaim

Let $\Cal G, \Cal I, \Cal A$, and $\Cal C$ be as in Lemma~1.3.  The
proof of Theorem~0.1 consists of an analysis of the pullback of these sets to
$E$. These pullbacks are denoted $G, J, A$, and $C$ respectively.  Claims
0,1, and 2 show how to homotop $E$ so that it efficiently intersects
$\Cal G, \Cal I, \Cal A$, and $\Cal C$.

\demo{Claim 0}  $E$ can be homotoped rel $\partial E$ to eliminate circle
components of
$\lambda\cap E$.  \enddemo

\demo{Proof of Claim 0}  Let $E_1,\cdots, E_n$ be a finite disjoint union of
embedded disks in $E$, which contain all the circle components of $\lambda|E$
and are disjoint from the non-circle leaves. The
$\pi_1$--injectivity of leaves of $\lambda$, implies that $E$ can be
homotoped rel $E-\cup E_i$ to eliminate the circle leaves, ie after
homotopy $E\cap\lambda=(E-\cup E_i)\cap \lambda$.  Notice that such a
homotopy does not increase the number of arc components of $E\cap \Cal
A$.   The resulting mapping of $E$ will be a branched immersion transverse
to $\lambda \cup \Cal A$.\qed\enddemo

\demo{Claim 1} $E$ can be homotoped rel $\partial E$ so that if $B$ is a
component of $E\cap \Cal A$, then $B$ is an arc which connects distinct leaves
of $\lambda|E$.  Furthermore $|B\cap C|\le 1$ and $\lambda|E$ still has
no circle components.\enddemo

\demo{Proof of Claim 1}  If $B$ is an innermost circle component of $E\cap
\Cal A$ let $E_1\subset \cirE$ be an embedded disk such that $B \subset
\cirE_1$ and $E_1 \cap (\lambda \cup \Cal A) =B$. The $\pi_1$--injectivity of
$\Cal A$ implies that a homotopy of $E$ supported in $E_1$ eliminates $B$ and
the resulting $E_1$ will be disjoint from $\lambda \cup A$. 

Now suppose that $B$ is an arc which connects the same leaves of
$\lambda|E$.  By passing to possibly another such one we can find an
embedded disk $E_1\subset E$ with $\partial E_1$ the union of $B$ and an
arc in a leaf of $\lambda\cap E$.  By Theorem~1.8 $B$ is a homotopically
inessential arc in $A$.  Since no leaf of $\lambda$ is isolated there is
small neighborhood $N(E_1)$ of $E_1$ such that $N(E_1)\cap A=B$.  Thus a
homotopy of $E$ supported in $N(E_1)$ eliminates this component $B$ of $E\cap
A$ without introducing other components of $E\cap A$.  Furthermore, the
homotopy can be carried out so that $\lambda|E$ continues to have no circle
components.  Therefore we can assume that if $B$ is a component of $E\cap
\Cal A$, then $B$ is an arc which connects distinct leaves of
$\lambda|E$.  Since $|E\cap C|$ is minimal it follows that $|B\cap C|\le 1$. 
\qed

\enddemo

\demo{Claim 2} After a homotopy of $E$ relative to $\partial E$, there is
no component $G_1$ of $G$ that is contained in $\cirE$ and is a disk
homeomorphic to $I\times I$ in such a way that $G_1\cap \partial_v G=\partial
I\times I$, where $0\times I$ is mapped to an  essential arc in $\Cal
A$.

\enddemo

\demo{Proof of Claim 2}  Suppose such a $G_1$ exists.  Let $K_i=i\times I$. 
By Theorem~1.8
$K_0$ cannot be homotoped rel $\partial K_0$ into a leaf of $\lambda$.   It
follows that $K_1$ is also essential.  An application of the Loop Theorem
shows that either $\Cal G$ contains an essential product disk or there exists
a relative homotopy in $\Cal G$ deforming $G_1$ into $\Cal A$.  (Ie a
homotopy $F\co (I \times I) \times I \to \Cal G$ such that $F_0 = f|G_1,\
F_t|i\times I=f|i\times I$ for $i\in \{0,1\},\ F_t|I\times i \subset
\partial_h|\Cal G$ for $i\in\{0,1\}$ and $F_1|I \times I \subset \Cal A$.) 
By Lemma~1.3 (3), the latter must occur.  Therefore after a homotopy of $E$
relative to $\partial E$, the number $|E\cap C|$ gets reduced by $2$
contradicting the minimality hypothesis. \qed

\enddemo

Figure~2.1 shows the $G_i$'s, in grey, and the $J_j$'s, in white, as subsets
$E$.  Points of $G \cap J$ which map to $\Cal C$ are indicated with dots. 
By Claim~2, regions like $G_1$ can be removed by an isotopy, and regions
like $G_2$ can not exist.  Note that components of $J|E$ need not be
$I$--bundles.  Indeed had $\lambda$ been a singular foliation $\mu$ split open,
then non $I$--bundle components of $J|E$ could have arisen from tangencies of
$D$ with $\mu$.

\hskip.95 truein
\relabelbox\small
\epsfxsize 3.3truein \epsfbox{Figure2.1}
\relabel {C}{$C$}
\relabel {J}{$J_j$}
\relabel {G1}{$G_1$}
\relabel {G2}{$G_2$}
\endrelabelbox
\vskip .05truein
\centerline{\small Figure 2.1}
\medskip

\hskip 1.5truein
\relabelbox
\epsfxsize 2.2truein \epsfbox{Figure2.2}
\relabel {G1}{$G_1$}
\relabel {G2}{$G_2$}
\relabel {G3}{$G_3$}
\relabel {C}{$C$}
\adjustrelabel <-2pt, 2pt > {F}{$F$}
\relabel {22}{$2$}
\endrelabelbox
\vskip 0truein
\centerline{\small Figure 2.2}
\medskip

If $F$ is a component of $E-\lambda$, define the index $I(F)$ to be
$1-\frac{1}{2}|F\cap \partial E|$.  There are only finitely many complementary
regions of non zero index by the Packet Lemma and Claim~1.  Since $E$ is a
disk,  $\sum_F I(F)=1$.  The regions of index $1/2$ are in 1--1 correspondence
with the outermost arcs of $\lambda\cap E$, moreover these regions are
components of $G$.  Let $N$ be the number of such regions.  Since
$1=\sum_F I(F)=\sum_{I(F)\le 0}I(F) + \sum_{I(F)>0}I(F)$ it follows that $N=
2 - 2\sum_{I(F)\le 0}I(F)$.
  
\demo{Claim 3}  $N>\frac{1}{3}|E\cap C| - \frac{1}{3\epsilon}\
\text{length}(\partial E)$.
\enddemo

\demo{Proof of Claim 3} Let $F$ be a closed complementary region of
$E-\lambda$ and let $G_1,\dots,G_p$ (respectively $J_1, \dots, J_q$) be
the components of
\ $G \cap F$ (respectively $J \cap F$).  See Figure~2.2.

Define $I(G_i) = 1 - \frac{1}{2}|G_i \cap (A \cup \partial E)|$ and
$I(J_i) = 1 - \frac{1}{2}|J_i \cap (A \cup \partial E)|$, thus $I(F) =
\sum_i I(G_i) + \sum_j I(J_j)$.  By construction $\partial E$ intersects $J$
only in interstitial fibres and by Theorem~1.8 such arcs cannot be homotoped,
fixing endpoints, into $\lambda$, thus

\roster
%\item If $I(F)> 0$ then $F\subset G$.

\item $I(J_j)\le 0$ for all j.

\item If $I(G_i)<0$ then $-I(G_i) \ge \frac{1}{6}|G_i \cap C|$.

\item If $I(G_i)=0$ then $|G_i \cap C|\le|G_i\cap\partial E|$.

\item If $I(G_i)>0$ then $G_i=F$, $I(G_i)= I(F)= \frac{1}{2}$ and $|G_i \cap
C| = 0$.
\endroster

(In the following sums we suppress double subscript notation by using the same
notation $\bigcup_i G_i \cup \bigcup_j J_j$ for the various different regions
$F$.) We have,
$$N= 2 - 2\sum_{I(F)\le 0}I(F)$$
$$>2\sum_{I(F)\le 0}\left(\sum_i(-I(G_i)) + \sum_j (-I(J_j))\right)$$
$$\ge 2\sum_{I(F)\le 0}\sum_i(-I(G_i))\tag{by 1}$$
$$= 2\sum_{I(F)\le 0}\sum_{I(G_i)<0}(-I(G_i))\tag{by 4}$$
$$\ge 2\sum_{I(F)\le 0}\sum_{I(G_i)<0}\frac{1}{6}|G_i \cap C|\tag{by 2}$$
$$= \frac{1}{3}\sum_{I(F)\le 0}\left(\sum_{I(G_i)\le 0} |G_i \cap C|-
\sum_{I(G_i)=0}|G_i \cap C|\right)$$
$$= \frac{1}{3}|E\cap C| - \frac{1}{3}\sum_{I(F)\le 0}\sum_{I(G_i)=0}
|G_i \cap C|\tag{by 4}$$
$$\ge\frac{1}{3}|E\cap C| - \frac{1}{3}\sum_{I(F)\le
0}\sum_{I(G_i)=0}|G_i \cap\partial E|\tag{by 3}$$
$$\ge\frac{1}{3}|E\cap C| - \frac{1}{3}\sum_F\sum_{G_i}|G_i
\cap\partial E|$$
$$=\frac{1}{3}|E\cap C| - \frac{1}{3}|G\cap\partial E|$$
$$\ge\frac{1}{3}|E \cap C| - \frac{1}{3\epsilon}\ \text{length}(\partial
E).\eqno{\hbox{(by Lemma~2.1)}\quad\qed}$$ 

\enddemo

\demo{Proof of Theorem~0.1}  

By the Ubiquity Theorem, \cite{G3}, there exists $ K,L>0$ such that if $D$ is
a least area disk such that $\text{length}(\partial D)/\text{area}(D)<L$,
then $$|E\cap \Cal C|/\text{area}(D)>K\tag{*}$$ 
where $E$ is any disk which spans $\partial D$.  This means that for
disks of small isoperimetric ratio, up to multiplicative constants, the 
wrapping number of $\partial D$ with $\Cal C$  is more or less the same as
$\text{area}(D)$.

Since there are $N$ regions of $E$ with index $\frac{1}{2}$, $\partial E \cap
G$ consists of at least $N$ components and therefore, 
$$\text{length}(\partial D)= \text{length}(\partial E) \ge \epsilon
N\tag{by Lemma~2.1}$$
$$\ge \frac{\epsilon}{3}|E\cap C| - \frac{1}{3}\text{length}(\partial
E)\tag{by Claim 3}$$
$$\ge \frac{K\epsilon}{3}\text{area}(D) - \frac{1}{3}\text{length}(\partial
D).\tag{by (*)}$$
It follows that 
$$\text{length}(\partial D) \ge \frac{K\epsilon}{12}\text{area}(D).\eqno\qed$$

\enddemo\eject

\enddocument